\documentclass{amsart}

\newtheorem{theorem}{Theorem}[section]

\theoremstyle{definition}
\newtheorem{definition}[theorem]{Definition}

\theoremstyle{remark}
\newtheorem{remark}[theorem]{Remark}

\numberwithin{equation}{section}

\newcommand{\abs}[1]{\lvert#1\rvert}

\begin{document}

\title{Nonexistence of global weak solutions for evolution equations with fractional Laplacian}

\author{Ahmad Z. Fino}
\address{Department of Mathematics, Faculty of Sciences, Lebanese University, P.O. Box 1352, Tripoli, Lebanon}
\email{ahmad.fino01@gmail.com; afino@ul.edu.lb}
\author{Evgeny I. Galakhov}
\address{RUDN University, 6 Miklukho-Maklaya street, Moscow, 117198, Russia}
\author{Olga A. Salieva}
\address{MSTU Stankin, 3a Vadkovsky Lane, Moscow, 127004, Russia}

\email{olga.a.salieva@gmail.com}







\begin{abstract}
In this paper, we are interested to analyze a nonlocal nonlinear parabolic equation with fractional Laplacian. We show that there are no nontrivial global weak solutions using the test function method. Nonexistence result of nontrivial weak solution for the corresponding elliptic case is obtained as a by-product.
\end{abstract}

\maketitle

\section{Introduction}
In this paper, we investigate the Cauchy problem for a nonlocal nonlinear
parabolic equation
\begin{equation}\label{pb1}
\left\{
\begin{array}{ll}
\,\,\displaystyle{u_t+(-\Delta)^{\beta/2}
(|u|^p)=|u|^{q}}&\displaystyle {x\in {\mathbb{R}^N},\;t>0,}\\
\\
\displaystyle{u(x,0)=u_0(x)}&\displaystyle{x\in {\mathbb{R}^N},}
\end{array}
\right. \end{equation} where $u_0\in L^1_{loc}(\mathbb{R}^N),$ $N\geq1,$
$0<\beta< 2,$ $p>0$, $q>1$, and the nonlocal operator
$(-\Delta)^{\beta/2}$ is defined by
\[(-\Delta)^{\beta/2}v(x):=\mathcal{F}^{-1}\left(\abs{\xi}^\beta
\mathcal{F}(v)(\xi)\right)(x)\]
for every $v\in
D((-\Delta)^{\beta/2})=H^\beta(\mathbb{R}^N),$ where
$H^\beta(\mathbb{R}^N)$ is the homogeneous Sobolev space of order
$\beta$ defined by
$$H^\beta(\mathbb{R}^N)=\left\{u\in L^2(\mathbb{R}^N);\;(-\Delta)^{\beta/2}u\in L^2(\mathbb{R}^N)\right\}$$
where $\mathcal{F}$ stands for the Fourier transform and
$\mathcal{F}^{-1}$ for its inverse, $\Gamma$ is
the Euler gamma function.

We start by the
\begin{definition}[Weak solution]${}$\\
Let $u_0\in L_{loc}^1(\mathbb{R}^N),$ $0<\beta<2$, $p>0$, $q>1$, and $T>0.$ We say that $u$ is a
weak solution of the problem (\ref{pb1}) if $u\in L^p((0,T),L^{2p}(\mathbb{R}^N))\cap L^q((0,T),L_{loc}^{q}(\mathbb{R}^N))$
and satisfies the equation
\begin{eqnarray}\label{weaksolution}
\int_{\mathbb{R}^N}u_0(x)\varphi(x,0)\,dx+\int_0^T\int_{\mathbb{R}^N}|u|^q\varphi(x,t)\,dx\,dt&=&\int_0^T\int_{\mathbb{R}^N}|u|^p(-\Delta)^{\beta/2}\varphi\,dx\,dt\nonumber\\
&-&\int_0^T\int_{\mathbb{R}^N}u\varphi_t(x,t)\,dx\,dt,
\end{eqnarray}
for all compactly supported $\varphi\in C([0,T],H^\beta(\mathbb{R}^N))\cap C^1([0,T],L^\infty(\mathbb{R}^N))$ such that $\varphi(\cdotp,T)=0$.
\end{definition}

Our main results are the following:

\begin{theorem}\label{theo1}
Let $0<\beta<2$, $p>0$, and $q>1$. For $u_0\in L^1_{loc}(\mathbb{R}^N)$, $u_0\geq0$, if
  $$p<q\leq p+\frac{\beta}{N},$$
then problem \eqref{pb1} has no nontrivial global weak solutions.
\end{theorem}

\begin{theorem}\label{theo2}
Let $0<\beta<2$, $p>0$, and $q>1$. For $u_0\in L^1_{loc}(\mathbb{R}^N),$ assume that there exists a constant $\varepsilon>0$ such that, for every $0<\gamma<N$, the initial datum verifies the following sign assumption:
$$u_0(x)\geq \varepsilon (1+|x|^2)^{-\gamma/2},$$
 If
 $$p<q< p+\frac{\beta}{\gamma},$$
 then problem \eqref{pb1} has no global weak solutions.
\end{theorem}

We prove these nonexistence theorems using a modification of the test function method of Zhang \cite{Zhang}, Mitidieri and Pohozaev \cite{PM1, P1}, and Kirane
et al. \cite{Finokirane,GK,KLT}; it was also used by Baras and Kersner \cite{Baras} for the study of the
necessary conditions for the local existence. However, it appears impossible to use this method in a direct way, since it was initially developed for other types of differential operators, such as integer powers of the Laplacian. But it is known \cite{CDV} that many problems containing these operators have relatively small solution sets. For instance, one cannot approximate with harmonic functions, in other words, with the solutions of the Laplace equation, an arbitrary function having interior maxima or minima (note that in the one-dimensional case harmonic functions are necessarily affine linear). This scarcity of solutions makes it easier to find additional nonlinear terms that exclude the existence of any nontrivial solutions.

In contrast, the sets of solutions for problems with fractional differential operators are normally much larger, sometimes even locally dense in $C(\mathbb{R}^n)$, as in the case of $s$-harmonic functions ($u$ such that $(-\Delta)^s u=0$), see \cite{DKK}, as well as in the case of higher order operators, see \cite{CDV, K}. Therefore, in order to obtain non-existence results in this situation one has to prevent the existence of any solutions from this larger set. Thus non-existence results in the fractional setting are always a delicate matter and require a substantial modification of the known techniques. As far as we know, up to now they were obtained only in some special cases, in particular in \cite{DKK} for systems of elliptic equations with fractional powers of the Laplacian, and by the second and third authors of the present paper in \cite{GS1, GS2, S} for some (mostly elliptic) inequalities and their systems with similar operators. Here we extend these results to the evolution problem (\ref{pb1}) not covered by previous results.

Noting that for the limit case $\beta\rightarrow 2$, Theorem \ref{theo1} is done in \cite[Theorem~29.1]{PM1}, while Theorem \ref{theo2} can be obtained by the same way.

Throughout this paper, positive constants will be denoted by $C$ and will change from line to line.

The rest of the paper consists of two sections and an appendix. In \S2 we prove Theorem \ref{theo1}, and in \S3, Theorem \ref{theo2}. In the appendix we give a proof of Ju's inequality used in \S2.


\section{Proof of Theorem \ref{theo1}}\label{section2}


The proof is by contradiction. Suppose that $u$ is a global weak solution to (\ref{pb1}), then, for all $T\gg1$, we have
\begin{eqnarray*}
\int_{\mathbb{R}^N}u_0(x)\varphi(x,0)\,dx+\int_0^T\int_{\mathbb{R}^N}|u|^q\varphi(x,t)\,dx\,dt&=&\int_0^T\int_{\mathbb{R}^N}|u|^p(-\Delta)^{\beta/2}\varphi\,dx\,dt\nonumber\\
&-&\int_0^T\int_{\mathbb{R}^N}u\varphi_t(x,t)\,dx\,dt,
\end{eqnarray*}
for all test function $\varphi\in C([0,T],H^\beta(\mathbb{R}^N))\cap C^1([0,T],L^\infty(\mathbb{R}^N))$ such that supp$\varphi$ is compact with $\varphi(\cdotp,T)=0$. \\
Now we take $\varphi(x,t):=\varphi_1^\ell(x) \varphi^\eta_2(t)$ with
$\varphi_1(x):=\Phi\left(\frac{|x|}{T^{\alpha}}\right),$ $\varphi_2(t):=\Phi\left(\frac{t}{T}\right),$
where
$\alpha=\frac{q-p}{\beta(q-1)}$, $\ell,\eta\gg 1$, and $\Phi$ a smooth nonnegative non-increasing function such that
\[
\Phi(r)=\left\{\begin {array}{ll}\displaystyle{1}&\displaystyle{\quad\text{if }0\leq r\leq 1/2,}\\\\
\displaystyle{0}&\displaystyle{\quad\text {if }r\geq 1.}
\end {array}\right.
\]
We have
\begin{eqnarray*}
\int_{\Omega_T}u_0(x)\varphi_1^\ell(x)\,dx+\int_0^T\int_{\Omega_T}|u|^q\varphi(x,t)\,dx\,dt&=&\int_0^T\int_{\mathbb{R}^N}|u|^p\varphi^\eta_2(t)(-\Delta)^{\beta/2}(\varphi_1^\ell(x))\,dx\,dt\nonumber\\
&-&\int_0^T\int_{\Omega_T}u \varphi_1^\ell(x)\frac{d}{dt}\varphi_2^\eta(t)\,dx\,dt,
\end{eqnarray*}
where
\begin{equation}\label{1a}
 \Omega_T:=\left\{ x\in{\mathbb{R}}^N;\;\;|x|\leq  T^{\alpha}\right\}.
\end{equation}
Using Ju's inequality $(-\Delta)^{\beta/2}
\left(\varphi_1^\ell\right)\leq
\ell\varphi_1^{\ell-1}(-\Delta)^{\beta/2}\left(\varphi_1\right)$
(see the Appendix), and the fact that $u_0\geq 0$, we obtain
\begin{eqnarray}\label{1}
&{}&\int_0^T\int_{\Omega_T}|u|^q\varphi(x,t)\,dx\,dt\nonumber\\
&{}&\leq C\int_0^T\int_{\Omega_T}|u|^p\varphi^\eta_2(t)\varphi_1^{\ell-1}(x)(-\Delta)^{\beta/2}(\varphi_1(x))\,dx\,dt\nonumber\\
&{}&+C \int_0^T\int_{\Omega_T}|u| \varphi_1^\ell(x)\varphi_2^{\eta-1}(t)\left|\frac{d}{dt}\varphi_2(t)\right|\,dx\,dt\nonumber\\
&{}&\leq C\int_0^T\int_{\Omega_T}|u|^p\varphi^{p/q}\varphi^{-p/q}\varphi^\eta_2(t)\varphi_1^{\ell-1}(x)\left|(-\Delta)^{\beta/2}(\varphi_1(x))\right|\,dx\,dt\nonumber\\
&{}&+C \int_0^T\int_{\Omega_T}|u| \varphi^{1/q}\varphi^{-1/q}\varphi_1^\ell(x)\varphi_2^{\eta-1}(t)\left|\frac{d}{dt}\varphi_2(t)\right|\,dx\,dt.
\end{eqnarray}
Therefore, using Young's inequality
\begin{equation}\label{2}
ab\leq\;\frac{1}{4}a^{q/p}\;+C\;
b^{\frac{q}{q-p}},
\end{equation}
with
\[
\left\{\begin{array}{l}
a=|u|^p\varphi^{p/q}\,\\
b=C \varphi^{-p/q}\varphi^\eta_2(t)\varphi_1^{\ell-1}(x)\left|(-\Delta)^{\beta/2}(\varphi_1(x))\right|\\
\end{array}
\right.
\]
in the first integral of the right-hand side of (\ref{1}), and the following Young's inequality
\begin{equation}\label{3}
ab\leq\;\frac{1}{4}a^{\;q}\;+C\;
b^{\;\tilde{q}}\qquad\text{where}\;\; \tilde{q}=\frac{q}{q-1},
\end{equation}
with
\[
\left\{\begin{array}{l}
a=|u| \varphi^{1/q},\\
b=C \varphi^{-1/q}\varphi_1^\ell(x)\varphi_2^{\eta-1}(t)\left|\frac{d}{dt}\varphi_2(t)\right|\\
\end{array}
\right.
\]
in the second integral of the right-hand side of (\ref{1}), we
get
\begin{eqnarray}\label{4}
&{}&\frac{1}{2}\int_0^T\int_{\Omega_T}|u|^q\varphi(x,t)\,dx\,dt\nonumber\\
&{}&\leq C\int_0^T\int_{\Omega_T}\varphi^\eta_2(t)\varphi_1^{\ell-\frac{q}{q-p}}(x)\left|(-\Delta)^{\beta/2}(\varphi_1(x))\right|^{\frac{q}{q-p}}\,dx\,dt\nonumber\\
&{}&\;+ C \int_0^T\int_{\Omega_T}\varphi_1^\ell(x)\varphi_2^{\eta-{\tilde{q}}}(t)\left|\frac{d}{dt}\varphi_2(t)\right|^{\tilde{q}}\,dx\,dt.
\end{eqnarray}
At this stage, we introduce the scaled variables: $s=
T^{-1}t,$ $y= T^{-\alpha}x$, in the right-hand side of (\ref{4}), we conclude that
\begin{equation}\label{5}
\int_0^T\int_{\Omega_T}|u|^q\varphi(x,t)\,dx\,dt \leq
\;C\;T^{-\delta},
\end{equation}
\noindent where $\delta=\frac{q}{q-1}-N\alpha-1=$. Now, noting that, as
$$
q\leq q^*:=p+\frac{\beta}{N}\quad\Longleftrightarrow\quad \delta\geq
0,
$$
\noindent we have to distinguish two cases:\\

\noindent $\bullet$ {\bf \underline{Case 1: $q<q^*$ (i.e. $\delta>0$)}}: By passing to the limit in (\ref{5}), as $T$ goes to $\infty$, we
get
\[
\lim_{T\rightarrow\infty}\int_0^{T}\int_{|x|\leq
T^{\alpha}}|u|^q\varphi(x,t)\,dx\,dt =0.
\]
\noindent Using the Lebesgue dominated convergence theorem and the fact that $\varphi(x,t)\rightarrow 1$ as $T\rightarrow\infty$, we conclude that
\[
\int_0^\infty\int_{\mathbb{R}^N}|u|^q(x,t)\,dx\,dt =0,
\]
and the by the continuity in time and space of $u$ we infer that $u\equiv0$.\\

\noindent $\bullet$ {\bf \underline{Case 2: $q=q^*$ (i.e. $\delta=0$) }}:  Using inequality
(\ref{5}) with $T\rightarrow\infty$ and taking into account the fact that $q=q^*,$ we have
$$
u\in L^q((0,\infty),L^q(\mathbb{R}^N));
$$
which implies that
\begin{eqnarray}\label{6}
&{}&\lim_{T\rightarrow\infty}\int_{T/2}^{T}\int_{|x|\leq
(BT)^{\alpha}}|u|^q\varphi(x,t)\,dx\,dt\nonumber\\
&{}& =\lim_{T\rightarrow\infty}\int_{0}^{T}\int_{|x|\leq
(BT)^{\alpha}}|u|^q\varphi(x,t)\,dx\,dt-\lim_{T\rightarrow\infty}\int_{0}^{T/2}\int_{|x|\leq
(BT)^{\alpha}}|u|^q\varphi(x,t)\,dx\,dt\nonumber\\
&{}&\quad=\int_{0}^{\infty}\int_{\mathbb{R}^N}|u|^q\varphi(x,t)\,dx\,dt-\int_{0}^{\infty}\int_{\mathbb{R}^N}|u|^q\varphi(x,t)\,dx\,dt=0.
\end{eqnarray}
On the other hand, repeating the same calculation as above by taking this time $\varphi_1(x):=
\Phi\left(\frac{|x|}{B^{\alpha}T^{\alpha}}\right),$ where
$1\leq B<T$ is large enough such that
when $T\rightarrow\infty$ we don't have $B\rightarrow\infty$ at the same time, and applying H\"older's inequality
\[
\int ab\leq\;\left(\int a^{\;q}\right)^{1/q}\;
\left(\int b^{\;\tilde{q}}\right)^{1/\tilde{q}},
\]
with
\[
\left\{\begin{array}{l}
a=|u| \varphi^{1/q},\\
b= \varphi^{-1/q}\varphi_1^\ell(x)\varphi_2^{\eta-1}(t)\left|\frac{d}{dt}\varphi_2(t)\right|\\
\end{array}
\right.
\]
in the second integral of the right-hand side of (\ref{1}) instead of Young's inequality, we arrive at
\begin{eqnarray*}
&{}&\int_0^T\int_{\Omega_{BT}}|u|^q\varphi(x,t)\,dx\,dt\\
&{}&\leq \frac{1}{4}\int_0^T\int_{\Omega_{BT}}|u|^q\varphi(x,t)\,dx\,dt+C\int_0^T\int_{\Omega_{BT}}\varphi^\eta_2(t)\varphi_1^{\ell-\frac{q}{q-p}}(x)\left|(-\Delta)^{\beta/2}(\varphi_1(x))\right|^{\frac{q}{q-p}}\,dx\,dt\\
&{}& +\; C\left(\int_{T/2}^T\int_{\Omega_{BT}}|u|^q\varphi(x,t)\,dx\,dt\right)^{1/q}\left( \int_0^T\int_{\Omega_{BT}}\varphi_1^\ell(x)\varphi_2^{\eta-{\tilde{q}}}(t)\left|\frac{d}{dt}\varphi_2(t)\right|^{\tilde{q}}\,dx\,dt\right)^{1/\tilde{q}},
\end{eqnarray*}
where
\[
\Omega_{BT}:=\left\{ x\in{\mathbb{R}}^N;\;|x|\leq  (BT)^{\alpha}\right\};
\]
therefore
\begin{eqnarray*}
&{}&\int_0^T\int_{\Omega_{BT}}|u|^q\varphi(x,t)\,dx\,dt\\
&{}&\leq C\int_0^T\int_{\Omega_{BT}}\varphi^\eta_2(t)\varphi_1^{\ell-\frac{q}{q-p}}(x)\left|(-\Delta)^{\beta/2}(\varphi_1(x))\right|^{\frac{q}{q-p}}\,dx\,dt\\
&{}& +\; C\left(\int_{T/2}^T\int_{\Omega_{BT}}|u|^q\varphi(x,t)\,dx\,dt\right)^{1/q}\left( \int_0^T\int_{\Omega_{BT}}\varphi_1^\ell(x)\varphi_2^{\eta-{\tilde{q}}}(t)\left|\frac{d}{dt}\varphi_2(t)\right|^{\tilde{q}}\,dx\,dt\right)^{1/\tilde{q}}.
\end{eqnarray*}
Thanks to the following rescaling:$\;\tau=
T^{-1}t,$ $\xi= (TB)^{-\alpha}x,$ taking into account the fact that $q=q^*$, we can easily conclude that
$$
\int_0^T\int_{\Omega_{BT}}|u|^q\varphi(x,t)\,dx\,dt \leq C\;B^{-1}+\; C\; \left(B^{N\alpha}\right)^{1/\tilde{q}}\left(\int_{T/2}^T\int_{\Omega_{BT}}|u|^q\varphi(x,t)\,dx\,dt\right)^{1/q}.
$$
Thus, taking the limits when $T\rightarrow\infty$ and then $B\rightarrow\infty,$ using (\ref{6}), we get:
$$
  \int_0^\infty\int_{\mathbb{R}^N}|u|^q(x,t)\,dx\,dt =0,
$$
i.e. $u\equiv0$. This completes the proof. $\hfill\square$\\

\section{Proof of Theorem \ref{theo2}}\label{section3}

The proof is also by contradiction by repeating the same argument in Section \ref{section2}. Suppose that $u$ is a global weak
solution to (\ref{pb1}), then,
\begin{eqnarray}\label{2.1}
&{}&\int_0^T\int_{\Omega_T}|u|^q\varphi(x,t)\,dx\,dt+\int_{\Omega_T}u_0(x)\varphi_1^\ell(x)\,dx\nonumber\\
&{}&\leq C\int_0^T\int_{\Omega_T}|u|^p\varphi^\eta_2(t)\varphi_1^{\ell-1}(x)(-\Delta)^{\beta/2}(\varphi_1(x))\,dx\,dt\nonumber\\
&{}&\quad+\,C \int_0^T\int_{\Omega_T}|u| \varphi_1^\ell(x)\varphi_2^{\eta-1}(t)\left|\frac{d}{dt}\varphi_2(t)\right|\,dx\,dt\nonumber\\
&{}&\leq C\int_0^T\int_{\Omega_T}|u|^p\varphi^{p/q}\varphi^{-p/q}\varphi^\eta_2(t)\varphi_1^{\ell-1}(x)\left|(-\Delta)^{\beta/2}(\varphi_1(x))\right|\,dx\,dt\nonumber\\
&{}&\quad+\,C \int_0^T\int_{\Omega_T}|u| \varphi^{1/q}\varphi^{-1/q}\varphi_1^\ell(x)\varphi_2^{\eta-1}(t)\left|\frac{d}{dt}\varphi_2(t)\right|\,dx\,dt.
\end{eqnarray}
Using the same test functions $\varphi(x,t)$ as above, we have
\begin{eqnarray*}
\int_{\Omega_T}u_0(x)\varphi_1^\ell(x)\,dx+\int_0^T\int_{\Omega_T}|u|^q\varphi(x,t)\,dx\,dt&=&\int_0^T\int_{\mathbb{R}^N}|u|^p\varphi^\eta_2(t)(-\Delta)^{\beta/2}(\varphi_1^\ell(x))\,dx\,dt\nonumber\\
&-&\int_0^T\int_{\Omega_T}u \varphi_1^\ell(x)\frac{d}{dt}\varphi_2^\eta(t)\,dx\,dt,
\end{eqnarray*}
where $\Omega_T$ is given by formula (\ref{1a}). Using Ju's inequality $(-\Delta)^{\beta/2}
\left(\varphi_1^\ell\right)\leq\ell\varphi_1^{\ell-1}(-\Delta)^{\beta/2}\left(\varphi_1\right)$, we get
\begin{eqnarray}\label{2.1a}
&{}&\int_0^T\int_{\Omega_T}|u|^q\varphi(x,t)\,dx\,dt+\int_{\Omega_T}u_0(x)\varphi_1^\ell(x)\,dx\nonumber\\
&{}&\leq C\int_0^T\int_{\Omega_T}|u|^p\varphi^\eta_2(t)\varphi_1^{\ell-1}(x)(-\Delta)^{\beta/2}(\varphi_1(x))\,dx\,dt\nonumber\\
&{}&+C \int_0^T\int_{\Omega_T}|u| \varphi_1^\ell(x)\varphi_2^{\eta-1}(t)\left|\frac{d}{dt}\varphi_2(t)\right|\,dx\,dt\nonumber\\
&{}&\leq C\int_0^T\int_{\Omega_T}|u|^p\varphi^{p/q}\varphi^{-p/q}\varphi^\eta_2(t)\varphi_1^{\ell-1}(x)\left|(-\Delta)^{\beta/2}(\varphi_1(x))\right|\,dx\,dt\nonumber\\
&{}&+C \int_0^T\int_{\Omega_T}|u| \varphi^{1/q}\varphi^{-1/q}\varphi_1^\ell(x)\varphi_2^{\eta-1}(t)\left|\frac{d}{dt}\varphi_2(t)\right|\,dx\,dt.
\end{eqnarray}
As $u_0(x)\geq \varepsilon (1+|x|^2)^{-\gamma/2}$, then
$$\int_{\Omega_T}u_0(x)\varphi_1^\ell(x)\,dx\geq \int_{|x|\leq T^\alpha/2}u_0(x)\,dx\geq\varepsilon\int_{|x|\leq T^\alpha/2} (1+|x|^2)^{-\gamma/2}\,dx\geq C \varepsilon T^{\alpha(N-\gamma)}.$$
Therefore, estimating the right-hand side of inequality (\ref{2.1a}) exactly in the same way as it was done with that of (\ref{1}) (see formulas (\ref{2})--(\ref{5})), we get
\begin{eqnarray}\label{2.4}
&{}&C \varepsilon T^{\alpha(N-\gamma)}+\frac{1}{2}\int_0^T\int_{\Omega_T}|u|^q\varphi(x,t)\,dx\,dt\nonumber\\
&{}&\leq C_1\int_0^T\int_{\Omega_T}\varphi^\eta_2(t)\varphi_1^{\ell-\frac{q}{q-p}}(x)\left|(-\Delta)^{\beta/2}(\varphi_1(x))\right|^{\frac{q}{q-p}}\,dx\,dt\nonumber\\
&{}&\quad+\, C_1 \int_0^T\int_{\Omega_T}\varphi_1^\ell(x)\varphi_2^{\eta-{\tilde{q}}}(t)\left|\frac{d}{dt}\varphi_2(t)\right|^{\tilde{q}}\,dx\,dt.
\end{eqnarray}
At this stage, we introduce again the scaled variables: $s=T^{-1}t,$ $y= T^{-\alpha}x$, in the right
hand-side of (\ref{2.4}), and, taking into account the non-negativity of the second term in the left-hand side of (\ref{2.4}), we conclude that
$$
C \varepsilon T^{\alpha(N-\gamma)} \leq
\;C_1\;T^{-\frac{q}{q-1}+N\alpha+1},
$$
that is
\begin{equation}\label{2.5}
\varepsilon \leq
\;C\;T^{-\delta^*},
\end{equation}
where $\delta^{*}=\frac{q}{q-1}-N\alpha-1+\alpha(N-\gamma)$. Now, noting that
$$
q< q^{**}:=p+\frac{\beta}{\gamma}\quad\Longleftrightarrow\quad \delta^{*}>0,
$$
then, by passing to the limit in (\ref{2.5}), as $T$ goes to $\infty$, we get a contradiction. $\hfill\square$\\

\begin{remark}
By applying the same calculation to the corresponding nonlocal elliptic equation
\begin{equation}\label{pb2}
(-\Delta)^{\beta/2}(|u|^p)=|u|^{q},\quad x\in \mathbb{R}^N,
 \end{equation}
 where $N\geq1,$
$0<\beta\leq 2,$ $p>0$, and $q>1$. We can obtain the following result: if
  $$p<q< \frac{Np}{(N-\beta)_+},$$
then problem (\ref{pb2}) has no nontrivial weak solutions.
\end{remark}


\section*{Appendix}
 In this appendix, we give a proof of Ju's inequality (see
Proposition $3.3$ in \cite{Ju}), in dimension $N\geq1$ where
$\delta\in[0,2]$ and $q\geq1,$ for all nonnegative Schwartz functions
$\psi$ (in the general case)
$$
(-\Delta)^{\delta/2}\psi^q\leq q\psi^{q-1}(-\Delta)^{\delta/2}\psi.
$$
The cases $\delta=0$ and $\delta=2$ are obvious, as well as $q=1.$ If $\delta\in(0,2)$ and $q>1,$
using \cite[Definition 3.2]{bogdan}, we have
$$
    (-\Delta)^{\delta/2}\psi(x)=-c_N(\delta)\int_{\mathbb{R}^N}
    \frac{\psi(x+z)-\psi(x)}{\abs{z}^{N+\delta}}\,dz,\quad
    \mbox{for all}\;x\in\mathbb{R}^N,
$$
where
$c_N(\delta)=2^{\delta}\Gamma((N+\delta)/2)/(\pi^{N/2}\Gamma(1-\delta/2)).$
Then
$$(\psi(x))^{q-1}(-\Delta)^{\delta/2}\psi(x)=-c_N(\delta)
\int_{\mathbb{R}^N}\frac{(\psi(x))^{q-1}\psi(x+z)-(\psi(x))^q}{\abs{z}^{N+\delta}}\,dz.$$
By Young's inequality we
have
$$(\psi(x))^{q-1}\psi(x+z)\leq\frac{q-1}{q}(\psi(x))^{q}+\frac{1}{q}(\psi(x+z))^q.$$
Therefore,
$$(\psi(x))^{q-1}(-\Delta)^{\delta/2}\psi(x)\geq
\frac{-c_N(\delta)}{q}\int_{\mathbb{R}^N}\frac{(\psi(x+z))^{q}-
(\psi(x))^q}{\abs{z}^{N+\delta}}\,dz=\frac{1}{q}(-\Delta)^{\delta/2}(\psi(x))^q.$$\\

\bibliographystyle{amsplain}

\end{document}